\documentclass[a4paper,11pt]{amsart}
\usepackage{amssymb}
\author[V.~S.~Luchko]{Victoriya S.~Luchko}
\address{
Algebra Department \\ Faculty of Mechanics and Mathematics\\ Kyiv
Taras Shevchenko University\\ Volodymyrskaia street, 01033 Kyiv,
Ukraine} \email{vsluchko@gmail.com}
\author[A.~P.~Petravchuk]{Anatoly P. Petravchuk}
\address{
Algebra Department \\ Faculty of Mechanics and Mathematics\\ Kyiv
Taras Shevchenko University\\ Volodymyrskaia street, 01033 Kyiv,
Ukraine} \email{aptr@univ.kiev.ua}
\title[On one-sided Lie nilpotent ideals of associative rings]
{On one-sided Lie nilpotent ideals of associative rings}
\date{March 7, 2008}

\DeclareMathOperator\differential{d}
\renewcommand\d\differential

\let\leq\leqslant
\let\geq\geqslant
\let\le\leq

\let\star *

\newtheorem{theorem}{Theorem}
\newtheorem{lemma}{Lemma}

\newtheorem{corollary}{Corollary}

\theoremstyle{definition}

\newtheorem{remark}{Remark}

\begin{document}

\sloppy

\begin{abstract}

  We prove that a
Lie nilpotent one-sided ideal of an associative ring $R$ is
contained in a Lie
 solvable two-sided ideal of $R$. An estimation of  derived length
 of such   Lie solvable ideal is obtained depending on the
 class of Lie nilpotency of the Lie nilpotent one-sided ideal of $R.$
  One-sided Lie nilpotent ideals
contained in ideals generated by commutators of the form $[\ldots
[ [r_1, r_{2}], \ldots ],  r_{n-1}], r_{n}]$ are also studied.

\end{abstract}

\maketitle

\section*{Introduction}

It is well-known that if $I$ is an one-sided nilpotent ideal of an
associative ring $R$ then $I$ is contained in a two-sided
nilpotent ideal of   $R$. Hence the following question is of
interest: for which one-sided ideal $I$ of the ring $R$ there
exists a two-sided ideal $J$ such that $J\supseteq I$ and $J$ has
properties like properties of $I$. In \cite{Petr} it was noted
that for an one-sided commutative ideal $I$ of a ring $R$ there
exists a nilpotent-by-commutative two-sided ideal $J$ of the ring
$R$ such that $J\supseteq I$.

Note  that Lie nilpotent and Lie solvable associative  rings were
 investigated by many authors (see, for example \cite{Jen}, \cite{Sharm},
 \cite{Streb}, \cite{Sys})
and the  structure of such rings  is studied well enough.

In this paper we prove that a Lie nilpotent one-sided ideal $I$ of
an associative ring $R$ is contained in a Lie solvable two-sided
ideal $J$ of $R$. An estimation (rather rough) of  Lie derived
length of the ideal  $J$ depending on  Lie nilpotency class  of
$I$ is also obtained (Theorem 1).

In  case when the Lie nilpotent one-sided ideal $I$ is contained
in the ideal $R_n$ of the ring $R$ generated by all commutators of
the form $[\ldots [ [r_1, r_{2}], \ldots ],  r_{n-1}], r_{n}]$ and
the Lie derived length
 of $I$ is  less then $n$ it is proved that
  $I$ is contained in a nilpotent two-sided ideal of $R$ (Theorem 2).

The notations in the paper are standard. If $S$ is a subset of an
associative ring $R$ then by $Ann^{l}_{R}(S)$ ($Ann^{r}_{R}(S)$)
we denote the left (respectively right) annihilator of $S$ in $R$.
We also denote by $R^{(-)}$  the  adjoint Lie ring of the
associative ring $R$. Further, by $R^{(-)}_{n}$  we denote the
$n$-th member of the lower central series of the Lie ring
$R^{(-)}$. Then $R_{n}=R^{(-)}_{n}+R^{(-)}_{n}\cdot R=$
$=R^{(-)}_{n}+R\cdot R^{(-)^{}}_{n}$ is a two-sided ideal of the
(associative) ring $R$. In particular, $R_2$ is a two-sided ideal
of the ring $R$ generated by all commutators of the form
$[r_{1},r_{2}]=r_{1}r_{2}-r_{2}r_{1}$, $r_{1},r_{2}\in R$. If $R$
is a  Lie solvable ring (i.e. such that $R^{(-)}$ is a solvable
Lie ring) then we denote  by $s(R)$  its  Lie derived length.
Analogously,  by $c(R)$ we denote   Lie nilpotency class  of a Lie
nilpotent ring  $R$.

\section{Lie nilpotent one-sided ideals}

\begin{lemma}\label{lemma:center} Let $I$ be an
one-sided ideal of an associative ring $R$ and $Z=Z(I)$ be the
center of  $I$. Then there exists an ideal $J$ in $R$ such that
$J^{2}=0$ and $[Z, R]\subseteq J$.
\end{lemma}
\begin{proof}
 Let, for example, $I$ be a right ideal from $R$. Take  arbitrary
elements $z\in Z$, $i\in I$, $r\in R$. Then it holds
$z(ir)-(ir)z=0$ (since $ir\in I$). This implies the equality
$i(zr-rz)=0$ since $z\in Z(I)$. As elements $z,i,r$ are
arbitrarily chosen then we have $I[Z,R]=0$. Consider the right
annihilator $T=Ann^{r}_{R}(I)$. It is clear that $T$ is a
two-sided ideal of the ring $R$ (since $I$ is  a right ideal of
$R$) what implies that $[Z,R]\subseteq T$.

Further,  for any element of the form $zr-rz$ from $[Z,R]$ and for
any $t\in T$ it holds $(zr-rz)t=z(rt)-r(zt)$. Since $rt\in T$ then
$z(rt)=0$. Besides,  $z\in I$ and therefore $zt=0$ what brings the
equality $(zr-rz)t=0$. It means that $[Z,R]\cdot T=0$.

Consider the left annihilator $J=Ann^{l}_{T}(T)$. It is easy to
see that $J$ is a two-sided ideal of the ring $R$. From relations
$[Z,R]\subseteq T$ and $[Z,R]\cdot T=0$ we have the inclusion
$[Z,R]\subseteq J$. It is also clear that $J^2=0$. Analogously one
can consider the case when $I$ is  a left ideal.
\end{proof}

\begin{theorem}\label{th:1}
 Let $R$ be an associative ring and $I$ be an
one-sided ideal of  $R$. If the subring $I$ is Lie nilpotent then
$I$ is contained in a Lie solvable two-sided ideal $J$ of  $R$
such that $s(J)\subseteq m(m+1)/2+m$ where $m=c(I)$ is  Lie
nilpotency class  of the subring $I$.
\end{theorem}
\begin{proof}
 Let for example $I$ be a right ideal. We prove our proposition by
the induction on the class of Lie nilpotency  $n=c(I)$ of the
subring $I$. If $n=1$ then $I$ is a commutative right ideal and by
 Lemma 1  the ring $R$ contains such  an ideal $T$ with
zero square that it holds $(I+T)/T\subseteq Z(R/T)$ in the
quotient ring $R/T$  where $Z(R/T)$ is the center of $R/T$. It
means that $I+T$ is a two-sided ideal of the ring $R$ and
$s(I+T)\leq 2$. Clearly $2=n+n(n+1)/2$ if $n=1$ and the statement
of Theorem is true in case $n=1$. Assume that the statement is
true in case $c(I)\leq n-1$ and prove it when $c(I)=n$. Denote by
$Z$ the center of the subring $I$. By Lemma 1 there exists an
ideal $T$ of $R$ with $T^2=0$ such that $[Z,R]\subseteq T$.
Consider the quotient ring $\overline{R}=R/T$. Then
$\overline{Z}=(Z+T)/T$ lies in the center of $\overline{R}$ and
therefore $\overline{Z}+\overline{Z}\cdot
\overline{R}=\overline{Z}+\overline{R}\cdot\overline{Z}$ is a
two-sided ideal of the ring $\overline{R}$. Since
$\overline{Z}\subseteq\overline{I}=(I+T)/T$  the ideal
$\overline{Z}+\overline{Z}\cdot\overline{R}$ is Lie nilpotent of
and its class of Lie nilpotency  $\le m.$ Further, the quotient
ring $\overline{R}/(\overline{Z}+\overline{Z}\cdot\overline{R})$
contains the right Lie nilpotent ideal
$\overline{I}+(\overline{Z}+\overline{Z}\cdot
\overline{R})/(\overline{Z}+\overline{Z}\cdot \overline{R})$ which
is Lie nilpotent of class of Lie nilpotency  $\leq m-1$. By the
induction assumption the last right ideal is contained in some Lie
solvable ideal of the ring
$\overline{R}/(\overline{Z}+\overline{Z}\cdot\overline{R})$ of
derived length $\leq\frac{(m-1)m}{2}+(m-1)$. Since
$\overline{Z}+\overline{Z}\cdot \overline{R}$ is Lie solvable and
its derived length  $\le m$ (even $\leq [log_{2}m]+1$ but we take
a rough estimation) and we consider  the quotient ring $R/T$ where
$T$ is Lie solvable of derived length 1, one can easily see  that
$I$ is contained in some Lie solvable (two-sided) ideal of derived
length which does not exceed
$$\frac{(m-1)m}{2}+(m-1)+(m+1)=\frac{(m+1)m}{2}+m.$$ Analogously
one can  consider the case when $I$ is right ideal.
\end{proof}

It seems to be unknown whether a sum of two Lie nilpotent
associative rings is Lie solvable. So the next statement can be of
interest (see also  results about sums of $PI$-rings in
\cite{Giam}).

\begin{corollary} Let $R$ be an associative ring which can decomposed into a
 sum $ R=A+B$ of its
Lie nilpotent subrings $A$ and $B$. If at least one of these
subrings is an one-sided ideal of $R$ then the ring $R$  is Lie
solvable.
\end{corollary}

\begin{remark} The statements of Theorem 1 and  its Corollary
become false when we replace Lie nilpotency of one-sided ideals by
Lie solvability. Really, consider full matrix ring
$R=M_{2}(\mathbb{K})$ over an arbitrary field $\mathbb{K}$ of
characteristic $\ne 2$. It is clear that
$$I=\left\{\left.\left(%
\begin{array}{cc}
  x & y \\
  0 & 0 \\
\end{array}%
\right)\right|x,y\in\mathbb{K}\right\}$$
 is a right Lie solvable
ideal of the ring $R$ but $I$ is  not contained in any Lie
solvable ideal of $R$ since $R$ is a non-solvable Lie  ring. It is
also clear that

$R=I+J$ where $J=\left\{\left.\left(%
\begin{array}{cc}
  0 & 0 \\
  z & t \\
\end{array}%
\right)\right|z,t\in\mathbb{K}\right\}$,

i.e. the simple associative  ring $R$ is a sum of two right Lie
solvable ideals.
\end{remark}

\section{On embedding
 of Lie nilpotent ideals in rings}

\begin{lemma}\label{lemma:power} Let $R$ be an associative ring,
$A$ be a Lie nilpotent subring of $R$ of Lie nilpotency class $<
m$. If $Z_0$ is a subring of $A$ such that $Z_{0}\subseteq Z(R)$
and $Z_{0}R\subseteq A$ then $Z_{0}^{m} R_{m}=0$.
\end{lemma}
\begin{proof}
Consider the two-sided ideal $J=Z_{0}+Z_{0}R=Z_{0}+RZ_{0}$ of the
ring $R$. As $J\subseteq A$ then $\underbrace{[J,...,J]}_{m}=0$ by
the condition $c(A)\leq m$. Further, it is easily to show that
$$[J,J]=[Z_{0}+Z_{0}R,Z_{0}+Z_{0}R]=Z^{2}_{0}[R,R].$$ By
induction on $k$ one can also  show that
$\underbrace{[J,...,J]}_{k}=Z^{k}_{0}\underbrace{[R,...,R]}_{k}$.
Then we have from the condition on $J$ that
$\underbrace{[J,...,J]}_{m}=Z^{m}_{0}\underbrace{[R,...,R]}_{m}=0$.
This implies the equality
$$Z^{m}_{0}R_{m}=Z^{m}_{0}(\underbrace{[R,...,R]}_{m}+\underbrace{[R,...,R]}_{m}\cdot
R)=\underbrace{[J,...,J]}_{m}+\underbrace{[J,...,J]}_{m}\cdot
R=0.$$
\end{proof}

\begin{lemma}\label{lemma:andr}Let $R$ be an associative ring,
$I$ be an ideal of  $R$. Then

1) if $J$ is a nilpotent ideal of the subring $I$ then $J$ lies in
a nilpotent ideal $J_I$ of the ring $R$ such that $J_{I}\subseteq
I$;

2) if $S=Ann^{l}_{I}(I)$ (or $Ann^{r}_{I}(I)$) then $S$ is
contained in a nilpotent ideal of the ring $R$ which is contained
in $I$.
\end{lemma}

The proof of this Lemma immediately follows from Lemma 1.1.5 from
\cite{Andr}.

\begin{theorem}\label{th:2} Let $R$ be an associative ring and
$I$ be a Lie nilpotent one-sided ideal of  $R$. If $I\subseteq
R_{n}$ and  Lie nilpotency class of $I$ is less then $n$ then $I$
is contained in an (associative)  nilpotent ideal of $R$.
\end{theorem}
\begin{proof}

Let for example $I$ be a right ideal of the ring $R$ and
$I\subseteq R_{n}$. One can assume that that $n\geq 2$ because the
statement of Theorem is obvious in case $n=1.$ We fix  $n\geq 2$
and prove the statement of Theorem by induction on the class of
Lie nilpotency $c=c(I)$ of the subring $I$. If $c=0$ then $I$ is
the zero ideal. and the proof is complete. Assume that the
statement is true for rings $R$ with $c(I)\leq c-1$ and prove it
in  case $c(I)=c$. Since $I$ is Lie nilpotent then by
Lemma~\ref{lemma:center} there exists a nilpotent ideal $T$ of the
ring $R$ such that in the quotient ring $\overline{R}=R/T$ it
holds $[\overline{Z_{0}},\overline{R}]=0$ where $Z_{0}$ is the
center of the subring $I$  and $\overline{Z_{0}}=(Z_{0}+T)/T$.
Then by Lemma~\ref{lemma:power} it holds the relation
$\overline{Z^{n}_{0}}\cdot \overline{R_{n}}=0.$
 If $\overline{Z^{n}_{0}}=0$
 then $\overline{Z_{0}}+\overline{Z_{0}}\overline{R}$ is a nilpotent
ideal of the ring $\overline{R}$ and then the subring $Z_{0}$ is
contained in the  nilpotent ideal $J=Z_{0}+T$ of the ring $R$.
Since in the quotient ring $R/J$ for the right ideal $(I+J)/J$ it
holds the inequality $c((I+J)/J)\leq c-1$ then by the inductive
assumption $(I+J)/J$ is contained in a nilpotent ideal $S/J$ of
the ring $R/J$. But then $I\subseteq S$ where $S$ is nilpotent
ideal of the ring $R$.

Let now  $\overline{Z^{n}_{0}}\neq 0$. Then
$\overline{Z^{n}_{0}}\subseteq
 Ann^{l}_{\overline{R_{n}}}(\overline{R_{n}})$ and since
 $\overline{Z_{0}}\subseteq \overline{R_{n}}$ then
 $\overline{Z^{n}_{0}}$ is
 contained in a nilpotent ideal $\overline{M}$ of the ring $\overline{R}$ by
 Lemma~\ref{lemma:andr}.
It is obvious that $\overline{Z_o}$ is nilpotent ideal of the ring
$\overline{R}$. Repeating the above considerations  we see that
$I\subseteq S$ where $S$ is a nilpotent ideal of the ring $R$.

\end{proof}

\begin{corollary}
 Let $R$ be an associative ring with condition  $R=[R,R]$. If $I$ is
 a Lie nilpotent one-sided ideal  of $R$ then there exists a nilpotent (two-sided) ideal $J$
  of the ring $R$ such that $I\subseteq J$

\end{corollary}

\begin{corollary}
 Let $R$ be a semiprime  ring. Then every Lie nilpotent one-sided
ideal is  contained in the center $Z(R)$ of the ring $R$ and has
trivial intersection with the ideal $R_2$.
\end{corollary}

\begin{proof}
Really since all nilpotent ideals of the ring $R$ are zero then by
Lemma~\ref{lemma:center} every Lie nilpotent one-sided ideal $I$
is contained in $Z(R)$. Since $IR\subseteq Z$ then
$[IR,R]=I[R,R]=0$.  Then from this equality we have $IR_{2}=I([R,
R]+[R, R]\cdot R)=0$. Denote $J=I\cap R_{2}$. It is easily to show
that $J\subseteq Ann^{l}_{R_{2}}(R_{2})$ and by
Lemma~\ref{lemma:andr} the intersection $J$ lies in a nilpotent
ideal of the ring $R$. As the ring $R$ is semiprime we have $J=0$.

\end{proof}

\end{document}